\documentclass{amsart}

\usepackage{geometry}
\usepackage{amssymb, amsmath}
\usepackage{bbm}
\usepackage{url}
\usepackage{enumitem}

\DeclareMathAlphabet{\dutchcal}{U}{dutchcal}{m}{n}

 \newtheorem{theorem}{Theorem}[section]
 \newtheorem{coro}[theorem]{Corollary}
 \newtheorem{lemma}[theorem]{Lemma}
 \newtheorem{proposition}[theorem]{Proposition}
 \theoremstyle{definition}
 \newtheorem{definition}[theorem]{Definition}
 \newtheorem{remark}[theorem]{Remark}

 \numberwithin{equation}{section}

\newcommand{\norm}[1]{\| #1 \|}
\newcommand{\bignorm}[1]{\bigl\| #1 \bigr\|}
\newcommand{\Bignorm}[1]{\Bigl\| #1 \Bigr\|}

\newcommand {\sfrac}[2] { {{}^{#1}\!\!/\!{}_{#2}}}
\newcommand {\pihalbe}{\sfrac{\pi}2}
\newcommand {\einhalb}{\sfrac{1}2}

\renewcommand{\Re}{\mathrm{Re}}
\renewcommand{\Im}{\mathrm{Im}}
\newcommand{\eins}{\mathbbm 1}

\newcommand{\EE}{\mathbb E}
\newcommand{\RR}{\mathbb R}
\newcommand{\CC}{\mathbb C}
\newcommand{\ZZ}{\mathbb Z}

\newcommand{\kernel}{\dutchcal{k}}
\newcommand{\BOUNDED}{\dutchcal{L}}
\newcommand{\DOMAIN}{\mathcal{D}}
\newcommand{\RANGE}{\mathcal{R}}
\newcommand{\calH}{\dutchcal{H}}
\newcommand{\calP}{\dutchcal{P}}
\newcommand{\calO}{\dutchcal{O}}
\newcommand{\calY}{\dutchcal{Y}}
\newcommand{\sector}[1]{S_{#1}}

\date{27 November 2026}

\title{The Hautus test for bounded analytic semigroups}

\author[Bernhard H. Haak]{Bernhard H. Haak}
\address{%
  Bernhard H. Haak \hfill\break\indent
  Institut de Math\'ematiques de Bordeaux\hfill\break\indent 
  Universit\'e de Bordeaux\hfill\break\indent 
  351, cours de la Lib\'eration\hfill\break\indent 
  33405 Talence cedex\hfill\break\indent 
  France}
\email{bernhard.haak@math.u-bordeaux.fr}

\author[Schwenninger]{Felix L.~Schwenninger}
\address{Felix L.~Schwenninger\hfill\break\indent 
Department of Applied Mathematics\hfill\break\indent
University of Twente \hfill\break\indent
P.O. Box 217 \hfill\break\indent
7500 AE Enschede, The Netherlands}
\email{f.l.schwenninger@utwente.nl}


\begin{document}

\begin{abstract}
  We prove that for bounded analytic semigroups on Hilbert spaces, the
  Hautus condition is sufficient for exact observability of an
  admissible observation operator, provided that the generator
  satisfies lower square function estimates. This complements a result
  of Le Merdy on the Weiss conjecture, where upper square function
  estimates play a corresponding role.  In general Banach spaces, we
  consider admissibility and exact observability in terms of norms of
  $\gamma$-radonifying operators and establish necessity and
  sufficiency of a randomised Hautus condition under suitable
  assumptions. For $L^p$-spaces, these results admit a formulation
  entirely in terms of Littlewood--Paley type inequalities.  As
  applications, we give an alternative proof of the characterisation
  of Riesz sequences of sectorial exponentials.
\end{abstract}

\maketitle

\section{Introduction}

This work is concerned with an observed Cauchy problem of the form
\[
  \left\{
  \begin{array}{lcl}
    x'(t) + A x(t) & = & 0 \\
    x(0) & = & x_0\\
    y(t) & = & C x(t)
  \end{array}\right.
\]
where ${-}A$ is the generator of a strongly continuous semigroup
$(T(t))_{t\ge 0}$ on a Banach space $X$ and $C:\DOMAIN(A)\to Y$ is an
observation operator into a Banach space $Y$, bounded with respect to
the graph norm of $A$.  More specifically, for a prescribed function
space $\calY$ for the output $t\mapsto y(t)$, we assume that $C$ is
$\calY$--admissible, which means that the observation map
$\Psi:x_0\mapsto y(\cdot)$ extends to a bounded operator $X\to\calY$,
where, in the Hilbert space setting, the canonical choice is
$\calY=L^2(\RR_+;Y)$.  The question of exact observability can then be
raised: is $\Psi$ bounded below, i.e., can one recover the initial
value from the observation in a stable manner?

\medskip\noindent In finite dimensions this observability question can
be answered by the well-known Kalman rank criterion.  Hautus, in
\cite{Hautus:test}, presented an alternative criterion: $C$ must be
injective on each eigenspace of $A$.  Unlike the Kalman rank
condition, Hautus' formulation lends itself to an extension to
infinite-dimensional systems, and such a criterion was proposed by
Russell and Weiss \cite{RussellWeiss:Hautus}. In our notation, their
criterion takes the form
\[
  \norm{x}^2 \lesssim
  \frac{\norm{(\lambda-A)x}^2}{(\Re\lambda)^2}
  +\frac{\norm{Cx}^2}{\Re\lambda},
  \qquad \lambda\in\CC_+,\quad x\in\DOMAIN(A),
\]
with a uniform constant. The authors showed the above criterion to be
necessary for exact observability and sufficient for approximate
observability of exponentially stable systems. The authors conjectured
that the criterion would also be sufficient for exact observability.
This conjecture was disproved a decade later by Jacob and Zwart
\cite{JacobZwart:counterexable-Weiss-conj}, who provided a
counterexample using an analytic diagonal semigroup on a suitably
chosen conditional basis.  The question has nevertheless continued to
attract considerable attention. Among the positive results, let us
mention the work of Jacob and Zwart
\cite{JacobZwart:Hautus-for-groups} for groups under certain growth
conditions. More recently, Daun and Jacob \cite{DaunJacob} gave a
positive answer for positive self-adjoint $A$ generating an
exponentially stable semigroup $(e^{-tA})_{t\ge0}$. Moreover, Daun and
Jacob found a counterexample for general (non-ananlytic) contraction
semigroups on Hilbert spaces.

\medskip\noindent In his work on the Weiss conjecture for bounded
analytic semigroups, Le Merdy \cite{LeMerdy:weiss-conj} proved that
the Weiss resolvent condition is sufficient for admissibility of an
observation operator $C$, provided that $A^{\einhalb}$ is itself an
admissible observation operator. In a complementary direction, one of
the main results of the present paper,
Theorem~\ref{thm:Hautus-Hilbert}, establishes the corresponding
statement for exact observability: for bounded analytic semigroups,
the Hautus condition is sufficient for an admissible observation
operator $C$ to be exactly observable, provided that $A^{\einhalb}$ is
itself exactly observable. In terms of square function estimates, the
two hypotheses on $A^{\einhalb}$ amount to upper and lower square
function estimates for $A$, respectively.

\medskip\noindent Outside the Hilbert space setting, the choice
$\calY=L^2(\RR_+;Y)$ is possible, though less compelling.  Experience
from harmonic analysis suggests that for $Y=L^p(\Omega)$, $p\not=2$,
Littlewood--Paley type norms
\[
  \Bignorm{    \Bigl(\int_0^\infty |CT(t)f|^2\,dt\Bigr)^\einhalb  }_{L^p}
\]
should be better suited to study the problem.  The appropriate
generalisation of these norms to arbitrary Banach spaces is provided
by norms of $\gamma$-radonifying operators. In earlier work, Haak and
Kunstmann have shown in \cite{HaakKunstmann:l-control-theory} that
Le\,Merdy's characterisation of admissibility has a perfect
counterpart in such $\gamma$-norms. For observability, these norms
prove useful again, so that we propose a randomised version of
the Hautus condition, for which we establish both necessity and
sufficiency under suitable assumptions in
Theorems~\ref{thm:gamma-Hautus-necess} and
\ref{thm:gamma-Hautus-suff}.  For $L^p$-spaces, these results can be
formulated entirely in terms of Littlewood--Paley type inequalities,
without any reference to Gaussian random variables, see
section~\ref{sec:LP-formulation}.

\section{Main result on Hilbert spaces}

We first recall some notions concerning sectorial operators and their
associated square functions.  For more background on the holomorphic
functional calculus for sectorial operators we refer to
\cite{Haase:Buch}. For $\theta\in(0,\pi)$ we write
\[
  \sector{\theta}=\{z\in\CC\setminus\{0\}: |\arg(z)|<\theta\}.
\]
A closed, densely defined operator $A$ on a Banach space $X$ is called
sectorial of angle $\omega$ if $\sigma(A)\subseteq\overline{S_\omega}$
and
\[
  \sup_{\lambda\notin\overline{S_\theta}} \norm{\lambda   R(\lambda,A)}<\infty
\]
for all $\theta\in(\omega,\pi)$.  Recall that $-A$ generates a bounded
analytic semigroup if and only if $A$ is sectorial of angle
$\omega<\pihalbe$. For $\theta>\omega$, let $H_0^\infty(S_\theta)$
denote the space of bounded holomorphic functions $\psi$ on $S_\theta$
satisfying
\[
  |\psi(z)|\lesssim \min\{|z|^\varepsilon,|z|^{-\varepsilon}\}
\]
for some $\varepsilon>0$.  For such functions, the Cauchy integral
\[
  \psi(A)=\frac1{2\pi i} \int_{\partial S_\nu} \psi(z)R(z,A)\,\mathrm{d} z,
  \qquad \omega<\nu<\theta,
\]
defines a bounded operator on $X$.

\medskip\noindent Throughout this paper, we write $\RR_+^*$ for the
positive half-line equipped with the multiplicative Haar measure
$\mathrm{d} t/t$, whereas $\RR_+$ and $\RR$ are equipped with Lebesgue
measure. If $X$ is a Hilbert space, the square function associated
with $\psi$ is given by
\[
  \norm{x}_{A,\psi}:= \Bigl(\int_0^\infty
  \norm{\psi(tA)x}_X^2\,\frac{\mathrm{d} t}{t} \Bigr)^\einhalb, \qquad x\in X,
\]
where we allow the value $+\infty$.  We say that $A$ admits upper
square function estimates if
\[
  \norm{x}_{A,\psi}\lesssim\norm{x}_X, \qquad x\in X,
\]
and lower square function estimates if
\[
  \norm{x}_X\lesssim\norm{x}_{A,\psi}, \qquad x\in X,
\]
for some non-zero $\psi\in H_0^\infty(S_\theta)$.  By a result of
McIntosh and Yagi \cite{McIntoshYagi:without}, for injective
sectorial operators on Hilbert spaces the square functions associated
with different non-zero functions in $H_0^\infty(S_\theta)$ are
equivalent.  In particular, the above estimates do not depend on the
particular choice of $\psi$. For a sectorial operator of angle
$\omega<\pihalbe$, a particularly natural choice is
$\psi(z)=z^\einhalb\,e^{-z}$, for which
\[
  \norm{x}_{A,\psi} = \Bigl(\int_0^\infty \norm{A^\einhalb   T(t)x}_X^2\,\mathrm{d} t \Bigr)^\einhalb.
\]
Thus, upper and lower square function estimates correspond precisely
to admissibility and exact observability of $A^\einhalb$,
respectively.  In the proof of our main result below, we will rather
work with functions of the form $\varphi_m(z)=z^m \, e^{-z}$.

\bigskip

In the following we prove that the Hautus test is sufficient for
bounded analytic semigroups on Hilbert spaces under the additional
assumption that $A$ admits lower square function estimates. This
extends recent works by Daun--Jacob \cite{DaunJacob} for self-adjoint
generators on Hilbert spaces and is in line with an earlier
counterexample by Jacob--Zwart
\cite{JacobZwart:counterexable-Weiss-conj} on the failure of the
Hautus test for exponentially stable analytic semigroups.

\begin{theorem}\label{thm:Hautus-Hilbert}{\bf (Hilbert space case)}
  Let $-A$ generate a bounded analytic semigroup $(T(t))_{t\ge0}$ on a
  Hilbert space $X$ and assume that $A$ satisfies lower square
  function estimates.  Then for any $L^{2}$-admissible
  $C\in\BOUNDED(\DOMAIN(A),Y)$ satisfying the \emph{Hautus condition},
\begin{equation}\label{eq:hautus}\tag{H}
  \norm{ x }^{2}  \lesssim \frac{1}{(\Re\lambda)^{2}}\norm{ (A-\lambda)x }^{2}+\frac{1}{\Re\lambda}\norm{ Cx }^{2}
\end{equation}
with uniform constants for all $\lambda\in\CC_{+}$ and all
$x\in \DOMAIN(A)$, it follows that $C$ is exactly observable.
\end{theorem}
\begin{proof}
  Let $\omega_A\in[0,\frac{\pi}{2})$ denote the sectoriality angle of
  $A$ and $\sector{\nu}$ be the sector
  $\{z\in\CC\colon |\arg(z)|\leq \nu\}$ for
  $\nu\in(\omega_A,\frac{\pi}{2})$.  We will show that there exists a
  suitable square function $\varphi\in H_0^\infty(\sector{\nu})$ such
  that
  \begin{equation}\label{eq:goal}
    \norm{\varphi(\cdot A)x}_{L^2(\RR_+^*;X)}     \lesssim    \norm{CT(\cdot)x}_{L^2(\RR_+;Y)},    \qquad x\in\DOMAIN(A),
  \end{equation}
  which readily implies the assertion by the assumed lower square
  function estimates. Here and in the following, $\RR_+^*$ denotes the
  positive half-line equipped with the multiplicative Haar measure
  $dt/t$, whereas $\RR_+$ and $\RR$ are equipped with Lebesgue
  measure.

  \medskip\noindent  \emph{Step 1: } 
  First observe that by sectoriality of $A$, the
  Hautus condition is equivalent to the inequality for $\lambda$ in a
  fixed sector $S{\nu}$.  By substituting $t=(\Re\lambda)^{-1}$ and
  $b=\Im(\lambda)t$ and putting $B=\tan\nu$, this means that
  \eqref{eq:hautus} is equivalent to
  \begin{equation}\label{eq:hautus2}
    c_B^2\norm{x}^2     \leq    \norm{(tA-1-ib)x}^2    +t\; \eins_{[-B,B]}(b)\norm{Cx}^2,
  \end{equation}
  for all $t>0$, $b\in\RR$ and $x\in\DOMAIN(A)$, where the constant
  $c_B>0$ depends only on $\nu$, the sectoriality constant of $A$ and
  the implicit constant in \eqref{eq:hautus}.

  \noindent The idea of the proof to use this form of the inequality
  for a suitably chosen element depending on $t$ and $b$ and integrate
  against $dt/t$ and $db$ such that the desired inequality
  \eqref{eq:goal} emerges.  More precisely, we aim to construct a
  family of operators $P_m(t,b)$ such that
  \begin{equation}\label{eq:ts1}
    \int_0^\infty\int_\RR     \norm{P_m(t,b)x}^2\,\mathrm{d}b\,    \frac{\mathrm{d}t}{t}
    \gtrsim
    \norm{\varphi_m(\cdot A)x}_{L^2(\RR_+^*;X)}^2,
  \end{equation}
  for a suitable square function $\varphi_m$, while the integral of
  $\norm{(tA-1-ib)P_m(t,b)x}^2$ against the same measure is
  comparatively small.  At the same time, the observation of
  $P_m(t,b)x$ should be controlled by that of $T(\cdot)x$.  For the
  following computations we first work with $x\in\DOMAIN(A)\cap R(A)$.
  Writing $x=Ay$ with $y\in\DOMAIN(A)$, analyticity of the semigroup
  gives, for $m\geq1$,
  \[
    \norm{(sA)^mT(s)x}    \lesssim    \min\{s\norm{Ax},s^{-1}\norm{y}\}.
  \]
  In particular, for $\varphi_m(z)=z^m\mathrm{e}^{-z}$, we have
  \[
    \varphi_m(\cdot A)x\in L^2(\RR_+^*;X).
  \]

  \medskip\noindent
  \emph{Step 2: }
  The requirements above motivate us to consider
  \[
    P_m(t,b)x
    =\int_0^\infty k_m(r,b)T(tr)x\,\mathrm{d}r,
    \qquad
    k_m(r,b)=\frac{\mathrm{d}^m}{\mathrm{d}r^m}
    \bigl(\kernel(r)\mathrm{e}^{ibr}\bigr),
  \]
  where $m\in\mathbb{N}$ and $\kernel(r)=r^{2m}\mathrm{e}^{-r}$.
  Integration by parts readily implies that
  \begin{equation}\label{eq:help1}
    P_m(t,b)x    =\int_0^\infty    \kernel(r)\mathrm{e}^{ibr}(tA)^mT(tr)x\,\mathrm{d}r.
  \end{equation}
  Note that, by analyticity of the semigroup and the vanishing of
  $\kernel$ and its derivatives at zero, all boundary terms
  vanish. Moreover, the latter integral converges in the graph norm of
  $A$.  Now Plancherel's theorem and Fubini's theorem give
  \begin{align}\notag
    \int_0^\infty\int_\RR
      \norm{P_m(t,b)x}^2\,\mathrm{d}b\,
      \frac{\mathrm{d}t}{t}
         = & \quad
      2\pi\int_0^\infty\int_0^\infty
      \kernel(r)^2\norm{(tA)^mT(tr)x}^2
      \,\mathrm{d}r\,\frac{\mathrm{d}t}{t}
      \notag\\
    = & \quad 
      2\pi\Bigl(\int_0^\infty
      \kernel(r)^2r^{-2m}\,\mathrm{d}r\Bigr)
      \Bigl(\int_0^\infty
      \norm{(sA)^mT(s)x}^2\,\frac{\mathrm{d}s}{s}\Bigr).
      \label{eq:second}
  \end{align}
  In particular, this yields \eqref{eq:ts1}.  Next, we consider
  $(tA-1-ib)P_m(t,b)x$, which by \eqref{eq:help1} and integration by
  parts rewrites to
  \begin{align*}
    (tA-1-ib)P_m(t,b)x
    =  &\;      (tA-1)\int_0^\infty    \kernel(r)\mathrm{e}^{ibr}(tA)^mT(tr)x\,\mathrm{d}r  \; + \;     \int_0^\infty    \frac{\mathrm{d}}{\mathrm{d}r}    \left(\kernel(r)(tA)^mT(tr)x\right)    \mathrm{e}^{ibr}\,\mathrm{d}r    \\
    = &  \;      \int_0^\infty    \bigl[(tA-1)\kernel(r)+\dot{\kernel}(r)    -(tA)\kernel(r)\bigr]    (tA)^mT(tr)x\,\mathrm{e}^{ibr}\,\mathrm{d}r    \\
    = &\;     \int_0^\infty    [\dot{\kernel}(r)-\kernel(r)]    \,\mathrm{e}^{ibr}(tA)^mT(tr)x\,\mathrm{d}r    =: E_m(t,b)x.
  \end{align*}
  By the same argument as above, using Plancherel's
  theorem, we get
\[
    \int_0^\infty\int_\RR
    \norm{E_m(t,b)x}^2\,\mathrm{d}b\,
    \frac{\mathrm{d}t}{t}
 \; = \;
    2\pi\Bigl(\int_0^\infty
    [\dot{\kernel}(r)-\kernel(r)]^2
    r^{-2m}\,\mathrm{d}r\Bigr)
    \Bigl(\int_0^\infty
    \norm{(sA)^mT(s)x}^2\,\frac{\mathrm{d}s}{s}\Bigr).
\]
  By the properties of the gamma function,
  \[
    \int_0^\infty r^j\mathrm{e}^{-\ell r}\,\mathrm{d}r
    =\frac{j!}{\ell^{j+1}},
  \]
  whence
  \[
    \int_0^\infty\kernel(r)^2r^{-2m}\,\mathrm{d}r
    =\frac{(2m)!}{2^{2m+1}}
    \quad\text{and}\quad
    \int_0^\infty
    [\dot{\kernel}(r)-\kernel(r)]^2r^{-2m}\,\mathrm{d}r
    =\frac{8m(2m-2)!}{2^{2m+1}}.
  \]
  Therefore,
  \begin{equation}\label{eq:Pm-Em-Hilbert}
    \norm{E_m(\cdot,\cdot)x}_{L^2(\RR_+^*\times\RR;X)}
    =
    \frac{2}{\sqrt{2m-1}}\,
    \norm{P_m(\cdot,\cdot)x}_{L^2(\RR_+^*\times\RR;X)}.
  \end{equation}

  \medskip\noindent
  \emph{Step 3: }
  We recall that
  \[
    P_m(t,b)x
    =\int_0^\infty k_m(r,b)T(tr)x\,\mathrm{d}r.
  \]
  By Minkowski's integral inequality and a change
  of variables,
  \begin{align}\notag
    \int_0^\infty\int_{-B}^B    \norm{CP_m(t,b)x}^2\,\mathrm{d}b\,\mathrm{d}t
    \leq & \;  \int_{-B}^B    \left(      \int_0^\infty |k_m(r,b)|      \left(\int_0^\infty
      \norm{CT(tr)x}^2\,\mathrm{d}t\right)^{\frac12}      \mathrm{d}r    \right)^2\mathrm{d}b    \notag\\
    = & \;     \norm{CT(\cdot)x}_{L^2(\RR_+;Y)}^2    \int_{-B}^B    \left(\int_0^\infty    |k_m(r,b)|r^{-\frac12}\,\mathrm{d}r\right)^2    \mathrm{d}b.
    \label{eq:first}
  \end{align}
  The latter integral is finite by the definition
  of $k_m$. Consequently,
  \begin{equation}\label{eq:observation-Pm-Hilbert}
    \norm{\eins_{[-B,B]}(b)\sqrt{t}\,CP_m(t,b)x}
    _{L^2(\RR_+^*\times\RR;Y)}
    \lesssim
    \norm{CT(\cdot)x}_{L^2(\RR_+;Y)}.
  \end{equation}

  \medskip\noindent
  \emph{Step 4: } 
  Applying \eqref{eq:hautus2} to $P_m(t,b)x$
  and integrating against $dt/t$ and $db$, we obtain
\[
    c_B^2
    \norm{P_m(\cdot,\cdot)x}_{L^2(\RR_+^*\times\RR;X)}^2
    \leq{}
    \norm{E_m(\cdot,\cdot)x}_{L^2(\RR_+^*\times\RR;X)}^2
    +
    \norm{\eins_{[-B,B]}(b)\sqrt{t}\,CP_m(t,b)x}_{L^2(\RR_+^*\times\RR;Y)}^2.
    \]
  Combining \eqref{eq:Pm-Em-Hilbert} and
  \eqref{eq:observation-Pm-Hilbert}, this yields
  \[
    \left(c_B^2-\frac{4}{2m-1}\right)
    \norm{P_m(\cdot,\cdot)x}_{L^2(\RR_+^*\times\RR;X)}^2
    \lesssim
    \norm{CT(\cdot)x}_{L^2(\RR_+;Y)}^2.
  \]
  Choosing $m$ sufficiently large such that
  $c_B^2-\frac{4}{2m-1}>0$, we conclude by
  \eqref{eq:second} that
  \[
    \norm{\varphi_m(\cdot A)x}_{L^2(\RR_+^*;X)}
    \lesssim
    \norm{P_m(\cdot,\cdot)x}_{L^2(\RR_+^*\times\RR;X)}
    \lesssim
    \norm{CT(\cdot)x}_{L^2(\RR_+;Y)}.
  \]
  Finally, the assumed lower square function estimates
  imply that $A$ is injective. By the mean ergodic
  decomposition of the bounded semigroup on the Hilbert
  space $X$, $A$ also has dense range.
  Hence $\DOMAIN(A)\cap R(A)$ is dense in $X$.
  By admissibility of $C$ and Fatou's lemma, the above
  estimate extends to all $x\in\DOMAIN(A)$, proving
  \eqref{eq:goal}. The assertion now follows from the
  assumed lower square function estimates.
\end{proof}

\begin{remark}
Note that the precise choice of the power $2m$ in
$\kernel(r)=r^{2m}\mathrm{e}^{-r}$ is crucial in the proof above. One can
show that the more general term $\kernel(r)=r^{c\,m}\mathrm{e}^{-r}$ with
$c\neq2$ does not imply that  the ratio
\[
  \frac{ \int_{0}^{\infty}[\dot{\kernel}(r)-\kernel(r)]^{2}r^{-2m}\,\mathrm{d}r }
       {\int_{0}^{\infty}\kernel(r)^{2}r^{-2m}\,\mathrm{d}r}
\]
tends to zero as $m\to\infty$.
\end{remark}

One might still ask wether the hypothsis of exact observability of $A^{\einhalb}$
is necessary for the conclusion in the Theorem. The answer is affirmative.

\begin{proposition}\label{prop:lower-sqf-necessary}
  Assume that ${-}A$ generates a bounded analytic semigroup
  $T(t)_{t\ge 0}$ on a Hilbert space $H$ and that the Hautus test
  \eqref{eq:hautus} implies exact observation for all
  $C \in \BOUNDED(D(A); Y)$. Then $A^\einhalb$ is exactly observable.
\end{proposition}
\begin{proof}
  For $\lambda\in \CC_+$, we have the elementary identity
  \[
    x = \frac1{\sqrt\lambda} (\sqrt\lambda+A^{1/2})^{-1}(\lambda-A)x    +\frac1{\sqrt\lambda}A^{1/2}x.
  \]
  But $A$ is sectorial of angle $\omega <\pihalbe$, and so
  $A^{\einhalb}$ is sectorial of angle $\omega/2<\sfrac{\pi}4$. Hence
  \[
     \norm{ (\sqrt\lambda+A^{1/2})^{-1} } \lesssim |\lambda|^{-1/2}
  \]
  with uniform constants for all $\lambda \in \CC_+$.  Consequently,
  \[
    \norm{ x } \lesssim
    \frac1{|\lambda|}\norm{ (\lambda-A)x }    +\frac1{\sqrt{|\lambda|}}\norm{ A^{1/2}x }.
  \]
  Since $ \Re\lambda\le|\lambda| $, this yields 
  \[
    \norm{ x }\lesssim \frac{\norm{ (\lambda-A)x     }}{\Re\lambda} +\frac{\norm{ A^{1/2}x     }}{\sqrt{\Re\lambda}},
  \]
  i.e. $A^{\einhalb}$ always satisfies the Hautus test and is hence
  exactly observable, i.e. $A$ admits lower square function estimates.
\end{proof}

The Proposition shows that the assumption of lower square function
estimates cannot be omitted from Theorem~\ref{thm:Hautus-Hilbert}. It
also highlights a feature of the counterexample by Jacob--Zwart
\cite{JacobZwart:counterexable-Weiss-conj} for exponentially stable
analytic semigroups: the operator $A$ in their counterexample does not
admit lower square function estimates.

\section{Main results in Banach spaces}

\subsection{Gaussian norms and square functions.}
Recall that throughout this paper, we write $\RR_+^*$ for the positive
half-line equipped with the multiplicative Haar measure
$\mathrm{d}t/t$, whereas $\RR_+$ and $\RR$ are equipped with Lebesgue
measure $\mathrm{d}t$. In this section we will work with the output
space $\calY=\gamma(\RR_+;Y)$.  We reproduce here a very brief summary
of the notions and notations we need below. For precise definitions,
background and further references on $\gamma$-functions and
$\gamma$-radonifying operators we refer to
\cite[Chap.~9]{FabulousFour:Band2}.

\smallskip \noindent Let $(S,\mu)$ be a $\sigma$-finite measure space
such that $L^2(S)$ is separable.  We say that a strongly measurable,
weakly square integrable function $f:S\to Y$ is a $\gamma$-function if
for every orthonormal basis $(e_n)$ in $L^2(S, \mu)$ the series
\[
  \sum_n\gamma_n\int_S f(t)e_n(t)\, \mathrm{d}\mu(t)
\]
converges in $L^1(\Omega;Y)$, and hence by Kahane's inequalities in
$L^p(\Omega; Y)$ for any $1\le p < \infty$. Here $(\gamma_n)$ is a
sequence of independent standard complex Gaussian random variables
over a probability space $\Omega$ and the integrals are understood in
the Pettis sense. We equip the space $\gamma(S;Y)$ of such functions
with the norm
\[
  \norm{f}_{\gamma(S;Y)}   :=  \EE\Bignorm{\sum_n\gamma_n   \int_S f(t)e_n(t)\,\mathrm{d}\mu(t)}_Y.
\]
By the rotational invariance of Gaussian sums, neither the definition
nor the norm depends on the choice of the orthonormal
basis. Equivalently, $f\in\gamma(S;Y)$ if the integral operator
\[
  R_f:L^2(S)\to Y,\qquad  R_f h=\int_S f(t)h(t)\,\mathrm{d}\mu(t),
\]
is well-defined and $\gamma$-radonifying. We will repeatedly use the
ideal property of the $\gamma$-norm: if $H_1,H_2$ are Hilbert spaces,
$V:H_1\to H_2$ and $U:X\to Y$ are bounded operators, and
$R\in\gamma(H_2;X)$, then
\[
  \norm{URV}_{\gamma(H_1;Y)} \leq
  \norm{U}\,\norm{R}_{\gamma(H_2;X)}\,\norm{V}.
\]
In particular, unitary transformations on the Hilbert space and
isometric embeddings of the Banach space preserve the
$\gamma$-norm. In particular, multiplication by a scalar function
$g\in L^\infty(S)$ defines a bounded operator on $\gamma(S;Y)$ of norm
at most $\norm{g}_\infty$. 

\medskip\noindent A family $\mathcal T\subseteq\BOUNDED(X;Y)$ is called $\gamma$-bounded
if there exists a constant $K$ such that
\[
  \EE\Bignorm{\sum_n\gamma_n T_nx_n}_Y   \leq K\EE\Bignorm{\sum_n\gamma_n x_n}_X
\]
for all finite choices of $T_n\in\mathcal T$ and $x_n\in X$. 

\medskip\noindent Let now $A$ be a densely defined $\gamma$-sectorial
operator of angle $\omega<\pi$ on a Banach space $X$.  By that we mean
that for all $\theta \in (\omega, \pi)$, the set of operators
\[
  \{ \lambda R(\lambda, A): \; \lambda \not \in \sector{\theta} \}
\]
is $\gamma$-bounded. In particular, $\gamma$-sectorial operators are
sectorial. The smallest such angle will be written by $\omega_\gamma(A)$.

\smallskip\noindent Let $\psi\in H_0^\infty(S_\theta)\setminus\{0\}$,
where $\theta>\omega_\gamma(A)$.  We say that $A$ admits upper
$\gamma$-square function estimates if
\[
  \norm{\psi(\cdot A)x}_{\gamma(\RR_+^*;X)} \lesssim_\psi \norm{x}_X,
  \qquad x\in\DOMAIN(A)\cap \RANGE(A),
\]
and lower $\gamma$-square function estimates if
\[
  \norm{x}_X\lesssim_\psi \norm{\psi(\cdot A)x}_{\gamma(\RR_+^*;X)},
  \qquad x\in\DOMAIN(A)\cap \RANGE(A).
\]
If $A$ has dense range, an upper square function estimate extends by
density to all $x\in X$.  By a result of Kalton and Weis
\cite[Proposition~7.7]{KaltonWeis:square-function-est}, the square
function norms associated with different non-zero functions
$\psi\in H_0^\infty(S_\theta)$ are equivalent on
$\DOMAIN(A)\cap R(A)$. As before, we will be particularly interested
in the case $\omega_\gamma(A)<\pihalbe$ and in functions of the form
$\varphi_\alpha(z)=z^\alpha\, e^{-z}$.  Notice that for Hilbert spaces
$X$, we have $\gamma(\RR_+^*; X) \simeq L^2(\RR_+^*; X)$ so that
$\gamma$-square function estimates fall back to those considered in
section \ref{sec:Hilbert}. The choice $\alpha=\einhalb$ connects
square-function estimates with admissibility and observability of
$A^\einhalb$ by simple change of measure:
\[
  \bignorm{ A^{\einhalb} T(\cdot)x }_{\gamma(\RR_+; X)} =
  \bignorm{ \varphi(tA) x }_{\gamma(\RR_+^*; X)}
\]


\bigskip \noindent For readers less familiar with $\gamma$-norms,
let us recall their close connection with classical Littlewood--Paley
theory. By the Khintchine--Maurey inequality, see e.g.
\cite[Theorem~7.2.13]{FabulousFour:Band2}, on any Banach lattice $Y$
of finite cotype we have
\[
  \EE\Bignorm{\sum_n\gamma_n y_n}_Y
  \simeq
  \Bignorm{\Bigl(\sum_n|y_n|^2\Bigr)^\einhalb}_Y.
\]
In particular, for $Y=L^p(\Omega)$, $1\le p<\infty$, the $\gamma$-norm
of the observation is equivalent to the familiar Littlewood--Paley
type norm
\[
  \norm{CT(\cdot)x}_{\gamma(L^2(\RR_+);L^p(\Omega))} \simeq \Bignorm{  \Bigl(\int_0^\infty |CT(t)x|^2\,\mathrm{d}t \Bigr)^\einhalb  }_{L^p(\Omega)}.
\]
More generally, by the ideal property of the $\gamma$-norm, one
obtains analogous identifications for Sobolev, Besov and
Triebel--Lizorkin spaces:
\[
  \begin{aligned}
    \gamma(L^2(S);W^{k,p})     &\simeq W^{k,p}(L^2(S)),\\
    \gamma(L^2(S);B^s_{p,q})   &\simeq B^s_{p,q}(L^2(S)),\\
    \gamma(L^2(S);F^s_{p,q})   &\simeq F^s_{p,q}(L^2(S)),
  \end{aligned}
\]
where $1\le p,q<\infty$. Thus, in all these examples, working with
$\gamma$-norms amounts to working with the corresponding
Hilbert-space-valued function spaces.

\smallskip \noindent Recall that a Banach space $Y$ is said to have
property $(\alpha)$ if
\[
  \EE\EE'\Bignorm{ \sum_{j,k} a_{jk}r_jr'_k y_{jk} }_Y \lesssim
  \EE\EE'\Bignorm{ \sum_{j,k} r_jr'_k y_{jk} }_Y
\]
uniformly for all finite families $(y_{jk})$ in $Y$ and scalars
$|a_{jk}|\leq 1$, where $(r_j)$ and $(r'_k)$ are independent sequences
of Rademacher variables. Banach lattices of finite cotype have this
property $(\alpha)$, see e.g. \cite[Corollary~7.6.20]{FabulousFour:Band2}.
Since property $(\alpha)$ is inherited by closed subspaces and
preserved under isomorphisms, it is shared by all the function spaces
mentioned above.

\subsection{Main results on $\gamma$--observability}

\begin{definition}
  We say that a bounded operator $C: \DOMAIN(A) \to Y$ is $\gamma$-admissible, if for all $x \in \DOMAIN(A)$,
  \[
    \bignorm{ CT(\cdot) x }_{\gamma(\RR_+; Y) } \;\lesssim \; \norm{x}_X
  \]
  Similarly, we say that $C$ is $\gamma$-observable if for all $x \in \DOMAIN(A)$,
  \[
    \bignorm{ CT(\cdot) x }_{\gamma(\RR_+; Y) } \;\gtrsim \; \norm{x}_X
  \]    
\end{definition}

\begin{theorem}[Haak~\cite{Haak:Diss}] \label{thm:gamma-Hautus-necess}
  {\bf (Necessity of the $\gamma$-Hautus test)}\hfill\\
  Let $X$ be a Banch space and $-A$ generator of a strongly continuous
  semigroup on $X$. Let $Y$ be a Banach space that has property
  $(\alpha)$.

  Assume that $C: \DOMAIN(A) \to Y$ is a $\gamma$-admissible operator that
  is also exactly $\gamma$-observable. Then there exists some $c>0$
  such that, for each $N \ge 1 $ and each choice of
  $\lambda_1, \ldots, \lambda_N$ with positive real part and each
  choice of $x_1, \ldots, x_N \in \DOMAIN(A)$,
  \begin{equation}    \label{eq:gamma-Hautus} \tag{$\gamma${-}H}
    c \, \EE \Bignorm{ \sum_{n=1}^N \gamma_n \Re(\lambda_n) x_n }_X
    \le
    \EE \Bignorm{ \sum_{n=1}^N \gamma_n (\lambda_n-A) x_n }_X
    +
    \EE \Bignorm{ \sum_{n=1}^N \gamma_n Re(\lambda_n)^\einhalb C  x_n }_Y
  \end{equation}
\end{theorem}

\noindent Before passing to the proof, we mention that due to Kahane's inequality,
\eqref{eq:gamma-Hautus} is equivalent to
\[
    c_p \, \Bigl( \EE \Bignorm{ \sum_{n=1}^N \gamma_n \Re(\lambda_n) x_n }_X^p \Bigr)^{\sfrac1p}
    \le
    \Bigl( \EE \Bignorm{ \sum_{n=1}^N \gamma_n (\lambda_n-A) x_n }_X^p  \Bigr)^{\sfrac1p}
    +
    \Bigl( \EE \Bignorm{ \sum_{n=1}^N \gamma_n Re(\lambda_n)^\einhalb C  x_n }_Y^p  \Bigr)^{\sfrac1p}
  \]
  where $0<p<\infty$. Moreover, letting $N=1$ shows that
  \eqref{eq:gamma-Hautus} implies the Hautus condition
  \eqref{eq:hautus}. Finally, for Hilbert spaces $X$, the
  choice $p=2$ and the orthogonality of the Gaussians implies that
  \eqref{eq:gamma-Hautus} is obtained by summing up the standard
  Hautus criterion \eqref{eq:hautus} $N$ times at different points, so
  that \eqref{eq:gamma-Hautus} and \eqref{eq:hautus} are indeed equivalent.

\begin{proof}
  By hypothesis, we have for all  $N \ge 1 $ and each choice of
  $\lambda_1, \ldots, \lambda_N$ with positive real part and each
  choice of $x_1, \ldots, x_N \in \DOMAIN(A)$,
  \begin{equation}    \label{eq:gamma-ex-obs-again}
    m\; \EE \Bignorm{ \sum_{n=1}^N  \gamma_n \Re(\lambda_n) x_n }_{\gamma(\RR_+; X)}
    \le
    m\; \EE \Bignorm{ \sum_{n=1}^N  \gamma_n \Re(\lambda_n) C T(\cdot) x_n }_{\gamma(\RR_+; Y)}.
  \end{equation}
  Using the fact that $\lambda-A$ generates the semigroup $(e^{\lambda t} T(t))_{t \ge 0}$ we see that
  \[
    T(t) x_0 = e^{-\lambda t } x_0 + \int_0^t e^{-\lambda (t-s)} T(s) (\lambda-A) x\,\mathrm{d}s.
  \]
  We use this equality for $x_0 \in \DOMAIN(A^2)$, and apply operator $C$ on
  both sides. This allows to split the right hand side of
  \eqref{eq:gamma-ex-obs-again} in two parts that we estimate
  separately.

  \medskip \noindent \emph{ First part: } let us write shorthand
  $e_z(t) := \eins_{\RR_+}(t) e^{tz}$. Then the operators
  $M_\lambda: \CC \to L^2(\RR_+)$ defined by
  $M_\lambda t := t \; \Re(\lambda)^\einhalb e_{-\lambda}$ are
  uniformly bounded over all $\lambda \in \CC_+$, each having norm
  one.  By The Haak-Kunstmann theorem, see
  \cite[Theorem~3.18]{HaakKunstmann:l-control-theory} or
  \cite[lookup]{FabulousFour:Band2}, the set of extenions
  $\{ M_\lambda: \lambda \in \CC_+\} \subset \BOUNDED(X; \gamma(\RR_+;
  X)$ is $\gamma$-bounded, with a constant that depends only on the
  property $(\alpha)$ constant. This immediately implies that
  \[
    \EE \Bignorm{ \sum_{n=1}^N \gamma_n \Re(\lambda_n) e_{-\lambda_n} \otimes C x_n }_{\gamma(\RR_+; Y)}
    \lesssim_{(\alpha)} 
    \EE \Bignorm{ \sum_{n=1}^N \gamma_n \Re(\lambda_n)^\einhalb  C x_n }_{Y}
  \]
  \medskip \noindent \emph{ Second part:} The convolution operators
  $K_\lambda f := f \ast \Re(\lambda) e_{-\lambda}$ are uniformly
  bounded on $L^2(\RR_+)$. By the Haak-Kunstmann theorem,
  \[
    \{ K_\lambda^\otimes : \Re(\lambda)>0 \} \subset \BOUNDED( \gamma(\RR_+; X) )
  \]
  are $\gamma$-bounded. But $K_\lambda^\otimes$ acts on
  $\gamma(\RR_+; X)$ by convolution, again, and so
  \begin{align*}
    & \; \EE \Bignorm{ \sum_{n=1}^N \gamma_n \bigl[ \Re(\lambda_n) e_{-\lambda_n}\bigr]  \ast  C T(\cdot) (\lambda_n-A) x_n }_{\gamma(\RR_+; Y)}\\
    \lesssim_{(\alpha)} & \; 
    \EE \Bignorm{ \sum_{n=1}^N \gamma_n C T(\cdot) (\lambda_n-A) x_n }_{\gamma(\RR_+; Y)}
    \le     M  \,   \EE \Bignorm{ \sum_{n=1}^N \gamma_n (\lambda_n-A) x_n }_{X}
  \end{align*}
  by the assumed $\gamma$-admissiblity with a constant, say $M$. Now
  we merge both parts, and obtain from \eqref{eq:gamma-ex-obs-again}
  \[
    \EE \Bignorm{ \sum_{n=1}^N  \gamma_n \Re(\lambda_n) C T(\cdot) x_n }_{\gamma(\RR_+; Y)}
    \le \EE \Bignorm{ \sum_{n=1}^N \gamma_n \Re(\lambda_n)^\einhalb  C x_n }_{Y}
    +  M \, \EE \Bignorm{ \sum_{n=1}^N \gamma_n (\lambda_n-A) x_n }_{X}
  \]
  whence the theorem.
\end{proof}

Before we turn to the new sufficiency result, let us mention that if
$A$ is $\gamma$-sectorial operator of angle
$\omega_\gamma < \pihalbe$, then $A^{\einhalb}$ does satisfy the
$\gamma$-Hautus test. This is shown along the lines of
Proposition~\ref{prop:lower-sqf-necessary}. Consequently, the
hypothesis of lower $\gamma$-square function estimates is necessary
for the following theorem to hold.

\begin{theorem}\label{thm:gamma-Hautus-suff}{\bf (Sufficiency of the $\gamma$-Hautus test)}\hfill\\
  Let $X, Y$ be Banach spaces and $A$ a densely defined
  $\gamma$-sectorial operator of angle $\omega_\gamma < \pihalbe$ with
  dense range that $A$ admits lower $\gamma$-square function
  estimates. Let $C: \DOMAIN(A)\to Y$ be a bounded $\gamma$-admissible
  observation operator. If the $\gamma$-Hautus condition
  \eqref{eq:gamma-Hautus} holds, the $C$ is exactly
  $\gamma$-observable.
\end{theorem}
\begin{proof}
  We follow the structure of the proof of Theorem~\ref{thm:Hautus-Hilbert}.
  While Step~1 requires substantial additional work in the
  Banach space setting, the remaining steps carry over by
  replacing the use of Plancherel's theorem and scalar
  $L^2$-estimates by the ideal property of the $\gamma$-norm.
  
  \medskip\noindent \emph{Step 1a: }
  If we choose in \eqref{eq:gamma-Hautus}
  \[
    \lambda_n=\frac{1+ib_n}{t_n},\qquad x_n=t_ny_n,
  \]
  where $y_n\in\DOMAIN(A)$ are arbitrary, then we obtain
  \[
      c\,\EE\Bignorm{\sum_{n=1}^N\gamma_n y_n}_X
      \le     \EE\Bignorm{        \sum_{n=1}^N\gamma_n(t_nA-1-ib_n)y_n}_X  +   \EE\Bignorm{    \sum_{n=1}^N\gamma_n\sqrt{t_n}\,Cy_n}_Y.
  \]
  Fix $\nu\in(\omega_\gamma,\pihalbe)$ and put $B=\tan\nu$.
  By $\gamma$-sectoriality, the resolvent  family 
  \[
    \dutchcal{R}_B    :=\bigl\{(tA-1-ib)^{-1}:t>0,\ |b|>B\bigr\}
  \]
  is $\gamma$-bounded. Indeed, with $\lambda=(1+ib)/t$,
  \[
    (tA-1-ib)^{-1}
    =-\frac{1}{1+ib}\lambda(\lambda-A)^{-1},
  \]
  and $\lambda$ lies outside the sector $K_\nu$ whenever
  $|b|>B$.
  Splitting the Gaussian sums into the indices with
  $|b_n|\le B$ and $|b_n|>B$, respectively, and using
  Gaussian contraction, \eqref{eq:gamma-Hautus} improves to
  \begin{equation}    \label{eq:almost-there}
      c_B\,\EE\Bignorm{\sum_{n=1}^N\gamma_n y_n}_X
      \le     \EE\Bignorm{   \sum_{n=1}^N\gamma_n(t_nA-1-ib_n)y_n}_X  +  \EE\Bignorm{ \sum_{n=1}^N\gamma_n  \eins_{[-B,B]}(b_n)\sqrt{t_n}\,Cy_n}_Y,
  \end{equation}
  where $c_B>0$ is independent of $N$ and of the chosen
  $(t_n,b_n)$ and $y_n$.

  \medskip\noindent \emph{Step 1b: }
  Let $Q\subset\RR_+\times\RR$ be a compact rectangle equipped with
  the measure $d\mu(t,b)=\frac{\mathrm{d}t}{t}\, db$.  Additionally to
  \eqref{eq:almost-there} obtained in Step~1, assume that we are given
  function $f:Q\to\DOMAIN(A)$ that is continuous with respect to the
  graph norm and assume that
  \[
    f,Af\in\gamma(L^2(Q,\mu);X).
  \]
  We claim that then a continuous $\gamma$-Hautus condition holds:
  \begin{equation}\label{eq:continuous-Hautus}
      c_B\bignorm{ f  }_{\gamma(L^2(Q,\mu);X)}    \le \;  \bignorm{ (tA-1-ib)f }_{\gamma(L^2(Q,\mu);X)} + \bignorm{ \eins_{[-B,B]}(b)\sqrt{t} \,Cf }_{\gamma(L^2(Q,\mu);Y)}.
  \end{equation}
  To prove the claim, we use a conditional expectation or averaging argument. Let
  $\calP=(Q_n)_{n=1}^N$ be a rectangular partition of $Q$ such
  that no rectangle $Q_n$ crosses either of the lines $b=-B$ or $b=B$.
  Put $w_n=\mu(Q_n)$ and choose points $(t_n,b_n)$ in the interiors of
  the $Q_n$.  Denote by $\Pi_{\calP}$ the orthogonal projection
  in $L^2(Q,\mu)$ onto the span of the orthonormal system
  $(w_n^{-\einhalb}\eins_{Q_n})_{n=1}^N$.  We use the same notation for the
  associated conditional averaging operator on vector-valued
  functions. Thus
  \[
    \Pi_{\calP}f=\sum_{n=1}^N\eins_{Q_n}f_n,      \qquad    f_n=\frac1{w_n}\int_{Q_n}f\,d\mu.
  \]
  We remark that $Cf \in\gamma(L^2(Q,\mu);Y)$, since
  $C f(t) = C(I+A)^{-1}\, (I+A) f(t)$.  Since the integrals may be
  taken in the graph norm of $A$,
  \[
    A\Pi_{\calP}f=\Pi_{\calP}(Af),     \qquad     C\Pi_{\calP}f=\Pi_{\calP}(Cf).
  \]
  Inserting $y_n=\sqrt{w_n}f_n$ into \eqref{eq:almost-there} and
  putting
  \[
    t_{\calP}=\sum_{n=1}^Nt_n\eins_{Q_n},    \qquad    b_{\calP}=\sum_{n=1}^Nb_n\eins_{Q_n},
  \]
  we obtain
  \begin{equation}
    \label{eq:jetzt-ist-erledigt}
      c_B\bignorm{ \Pi_{\calP}f }_\gamma    \le \;    \bignorm{ t_{\calP}\Pi_{\calP}(Af)    -(1+ib_{\calP})\Pi_{\calP}f }_\gamma
          \;+ \;         \bignorm{ \eins_{[-B,B]}(b)\sqrt{t_{\calP}}\,         \Pi_{\calP}(Cf) }_\gamma.
  \end{equation}
   Here we have used that $\eins_{[-B,B]}$ is constant on each cell,
   up to null sets, and that the $\gamma$-norm actually calculates
   with the finite orthonormal system
   $(w_n^{-\einhalb}\eins_{Q_n})_{n=1}^N$.  As the mesh of the partition
   tends to zero, $\Pi_{\calP}$ converges strongly to the identity on
   $L^2(Q,\mu)$. For any Banach space $E$ and suitable $g: Q \to E$,
   denote by $R_g$ the associated $\gamma$-radonifying operator. Since
   $\Pi_{\calP}$ is an orthogonal projection,
  \[
    (R_g\Pi_{\calP})^*x^*     =\Pi_{\calP}R_g^*x^*     \longrightarrow R_g^*x^*
  \]
  for every $x^*\in E^*$, and
  \[
    \bignorm{(R_g\Pi_{\calP})^*x^*}_{L^2(Q,\mu)}     \le \bignorm{R_g^*x^*}_{L^2(Q,\mu)}.
  \]
  Hence the $\gamma$-dominated convergence theorem
  \cite[Theorem~9.4.2]{FabulousFour:Band2} yields
  \[
    \Pi_{\calP}g\longrightarrow g \quad\text{in }\gamma(L^2(Q,\mu);E).
  \]
  Moreover, $t_{\calP}$, $b_{\calP}$ and $\sqrt{t_{\calP}}$ converge
  uniformly to $t$, $b$ and $\sqrt t$, respectively. Since scalar
  $L^\infty$-multipliers act boundedly on spaces of
  $\gamma$-radonifying operators, the claim follows upon passing
  to the limit in \eqref{eq:jetzt-ist-erledigt}.

  \medskip\noindent \emph{Step 1c: }
  For $x\in\DOMAIN(A)\cap \RANGE(A)$ and
  $\varphi_m(z)=z^m\, e^{-z}$ we verify that the associated square
  functions are indeed $\gamma(\RR_+^*; X)$-radonifying. Writing
  $x=Ay$ with $y\in\DOMAIN(A)$, analyticity of the semigroup gives
  \[
    \norm{ \varphi_m(sA)x }_X  \lesssim   \min\{s \norm{Ax}_X, s^{-1} \norm{y}_X\}.
  \]
  Since
  \[
    s\frac{\mathrm{d}}{\mathrm{d}s}\varphi_m(sA)x   =m\varphi_m(sA)x-\varphi_{m+1}(sA)x,
  \]
  putting $g(u) := \varphi_m(e^u A) x$, $g \in W^{1,1}(\RR;X)$ and
  hence, by Corollary~\ref{coro:W11}
  \[
    t \mapsto \varphi_m(t A)x   \in \gamma( L^2 (\RR_+^*); X ).
  \]

  \medskip\noindent\emph{Step 2:}
  We recall some notation from the proof of  Theorem~\ref{thm:Hautus-Hilbert}.
  \begin{equation}
    \label{eq:defPm-and-Em}
    \begin{aligned}
    P_m(t, b) x :=  & \;  \int_0^\infty  \kernel(r) \, e^{ibr}(tA)^m T(tr) x\,\mathrm{d}r   \\
    E_m(t, b) x :=  &\;   \int_{0}^{\infty}[\dot{\kernel}(r)-\kernel(r)]\,\mathrm{e}^{ibr}\,(tA)^{m}T(tr)z\,\mathrm{d}r
    \end{aligned}
  \end{equation}
  where $\kernel(r) = r^{2m}e^{-r}$.  Also recall that, by integration
  by parts,
  \begin{equation}    \label{eq:Em-Pm}
        E_m(t,b)x=(tA-1-ib)P_m(t,b)x.
      \end{equation}
      We consider the Hilbert space
  $\calH := L^2(\RR_+^*\times\RR)$ and extend the kernel $\kernel$ by
  zero for $r<0$.  Fourier transformation in the second variable and
  the subsequent change of variables $s=tr$ correspond to a unitary
  transformation on $\calH$.  By the ideal property of the
  $\gamma$-norm, and observing that $P_m$ is an inverse Fourier
  transform (see \eqref{eq:defPm-and-Em}), we therefore obtain
  \begin{equation}    \label{eq:gamma-norm-on-Pm}
    \begin{aligned}
      \bignorm{ P_m(\cdot,\cdot)x }_{\gamma(\calH;X)}    &=\sqrt{2\pi}\,       \bignorm{ (t, r) \mapsto   \kernel(r)(tA)^mT(tr)x }_{\gamma(\calH;X)}\\
      &=\sqrt{2\pi}\,  \bignorm{ (s, r) \mapsto \kernel(r)r^{-m}\varphi_m(sA)x }_{\gamma(\calH;X)}\\
      &=\sqrt{2\pi}\,    \bignorm{ \kernel(r)r^{-m} } _{L^2(\RR_+)}   \norm{ \varphi_m(\cdot A)x  }_{\gamma(\RR_+^*;X)}.
    \end{aligned}
  \end{equation}
  Similarly,
  \[
      \bignorm{ E_m(\cdot,\cdot)x }_{\gamma(\calH;X)}  =\sqrt{2\pi}\,  \bignorm{ (\dot{\kernel}(r)-\kernel(r))r^{-m}    }_{L^2(\RR_+)}
      \bignorm{ \varphi_m(\cdot A)x }_{\gamma(\RR_+^*;X)}.
  \]
  In particular, the computation in the the proof of the Hilbert space case in Theorem~\ref{thm:Hautus-Hilbert}
  yields
  \begin{equation}\label{eq:gamma-error-small}
    \bignorm{ E_m(\cdot,\cdot)x }_{\gamma(\calH;X)}   = \frac{2}{\sqrt{2m-1}} \bignorm{ P_m(\cdot,\cdot)x }_{\gamma(\calH;X)}.
  \end{equation}

  \medskip\noindent \emph{Step 3: }
  Recall that
  \[
    P_m(t,b)x=\int_0^\infty k_m(r,b) \, T(tr)x\,\mathrm{d}r,
    \quad\text{where}\quad
    k_m(r,b)=\frac{d^m}{d r^m} \bigl(\kernel(r)e^{ibr}\bigr),
  \]
   and consider the scalar integral operator
  \[
    J_m:L^2(\RR_+)\longrightarrow \calH,\qquad    (J_mh)(t,b)=\eins_{[-B,B]}(b) \;     \int_0^\infty k_m(r,b)\, \sqrt{t} \,h(tr)\,\mathrm{d}r.
  \]
  By Minkowski's inequality and a change of variables,
  \[
    \begin{aligned}
      \norm{ J_mh }_\calH   &\le \int_0^\infty \norm{ k_m(r,\cdot) }_{L^2(-B,B)} \norm{\sqrt t\,h(tr) }_{L^2(\frac{\mathrm{d}t}{t})}      \,\mathrm{d}r\\
      &=\norm{ h  }_{L^2(\RR_+)}      \int_0^\infty   \norm{ k_m(r,\cdot) }_{L^2(-B,B)}  \frac{\mathrm{d}r}{\sqrt r}.
    \end{aligned}
  \]
  The latter integral is finite. Since $J_m$ induces a bounded
  operator from $\gamma(L^2(\RR_+);Y)$ into $\gamma(H;Y)$, we have
  \begin{equation}\label{eq:gamma-observation-term}
    \norm{   \eins_{[-B,B]}(b)\sqrt t\,CP_m(t,b)x }_{\gamma(H;Y)}     \le \norm{ J_m } \, \norm{ CT(\cdot)x }_{\gamma(L^2(\RR_+);Y)}.
  \end{equation}

  \medskip\noindent \emph{Step 4: }
  On every compact rectangle $Q\subset\RR_+\times\RR$, the function
  $P_m(\cdot,\cdot)x$ is continuous with respect to the graph
  norm. Moreover, using \eqref{eq:Em-Pm}
  \[
    AP_m(t,b)x=\frac1t   \bigl(E_m(t,b)x+(1+ib)P_m(t,b)x\bigr)
  \]
  which implies that $AP_m(\cdot,\cdot)x$ is $\gamma$-radonifying on
  $Q$.  We may thus apply \eqref{eq:continuous-Hautus} to
  $f(t,b)=P_m(t,b)x$ on each compact rectangle
  $Q_N=[N^{-1},N]\times[-N,N]$.  Since all three functions involved
  are globally $\gamma$-radonifyin the $\gamma$-dominated convergence
  theorem \cite[Corollary 9.4.3]{FabulousFour:Band2}, allows us to
  pass to the limit $N\to\infty$.  We obtain
  \[
             c_B \bignorm{ P_m(\cdot,\cdot)x }_{\gamma(\calH;X)}
    \le
             \bignorm{ E_m(\cdot,\cdot)x }_{\gamma(\calH;X)}
      \;+\;  \bignorm{ \eins_{[-B,B]}(b)\sqrt t\,CP_m(t,b)x }_{\gamma(\calH;Y)}.
  \]
  Combining \eqref{eq:gamma-error-small} and
  \eqref{eq:gamma-observation-term}, this yields
  \[
    \left(c_B-\frac{2}{\sqrt{2m-1}}\right)     \norm{ P_m(\cdot,\cdot)x }_{\gamma(\calH;X)}
    \le \norm{ J_m }\,   \norm{ CT(\cdot)x }_{\gamma(L^2(\RR_+);Y)}.
  \]
  Choosing $m$ sufficiently large, and using the assumed
  lower $\gamma$-square function estimates, we conclude by \eqref{eq:gamma-norm-on-Pm}
  that
  \[
    \norm{ x }_X    \lesssim     \norm{ \varphi_m(\cdot A)x    }_{\gamma(L^2(\frac{\mathrm{d}s}{s});X)}
    \lesssim
     \bignorm{ P_m(\cdot,\cdot)x }_{\gamma(\calH;X)} 
    \lesssim
    \norm{ CT(\cdot)x }_{\gamma(L^2(\RR_+);Y)}. \qedhere
  \]
\end{proof}

\subsection{ A Littlewood–Paley formulation of the main theorem}\label{sec:LP-formulation}\hfill\\
It may be useful to write out the previous result in the spacial case
of $L^p$-spaces, since then we can suppress all Gaussian sums in
favour of vector-valued inequalities à la Garcia-Cuerva and Rubio De
Francia. Recall that, for $1\le p<\infty$, Kahane–Khintchine and
Fubini yield
\[
  \EE \Bignorm{ \sum_n\gamma_n f_n }_{L^p} \asymp_p \Bignorm{ \Bigl(\sum_n|f_n|^2\Bigr)^{\einhalb}  }_{L^p}.
\]
The continuous counterpart is the identification
$\gamma(L^2(S);L^p(\Omega)) \simeq L^p(\Omega;L^2(S)).$ In particular,
this includes $p=1$. Moreover, all these $L^p$-spaces have property
$(\alpha)$, so both the necessary and the sufficient parts of
Theorems~\ref{thm:gamma-Hautus-necess} and \ref{thm:gamma-Hautus-suff}
are available.  Let $1\le r,p<\infty$, let $X=L^r(\Omega)$ and
$Y=L^p(\calO)$, and let $-A$ be the generator of a bounded analytic
semigroup $T(t)_{t\ge 0}$ that has dense range and satisfes 
\[
\Bignorm{ \Bigl(\sum_n|z_n(z_n-A)^{-1}x_n|^2\Bigr)^{1/2} }_r
\lesssim_\nu \Bignorm{ \Bigl(\sum_n|x_n|^2\Bigr)^{1/2}}_r
\]
for $|\arg(z_n)|>\theta > \omega_\gamma(A)$ where
$\omega_\gamma(A)<\pihalbe$.  Assume that $A$ satisfies the lower Littlewood-Paley 
square-function estimate
\[
  \norm{ x }_{L^r(\Omega)} \lesssim \Bignorm{ \Bigl(\int_0^\infty
    |A^{\einhalb}T(t)x|^2\,dt\Bigr)^{\einhalb} }_{L^r(\Omega)}.
\]
Let $C:D(A)\to L^p(\calO)$ be a bounded, Littlewood-Paley admissible observation in the sense that
\[
  \Bignorm{ \Bigl(\int_0^\infty |CT(t)x|^2\,dt\Bigr)^{\einhalb}
  }_{L^p(\calO)} \lesssim\norm{ x }_{L^r(\Omega)}. 
\]
Then the following assertions are equivalent.
\begin{enumerate}
\item  Exact observability in the Littlewood–Paley sense:
  \[
    \norm{ x }_{L^r(\Omega)} \lesssim \Bignorm{  \Bigl(\int_0^\infty |CT(t)x|^2\,dt\Bigr)^{\einhalb}  }_{L^p(\calO)}.  
  \]
\item The vector-valued Hautus condition: there exists $c>0$ such that
  for every finite sequence $(x_n)_n\subset D(A)$ and every finite
  sequence $(\lambda_n)_n\subset\CC_+$,
  \[
    \Bignorm{  \Bigl(\sum_n(\Re\lambda_n)^2|x_n|^2\Bigr)^{\einhalb}  }_{L^r(\Omega)}
    \lesssim \Bignorm{  \Bigl(\sum_n|(\lambda_n-A)x_n|^2\Bigr)^{\einhalb}  }_{L^r(\Omega)}
    \;+\; \Bignorm{  \Bigl(\sum_n\Re\lambda_n\,|Cx_n|^2\Bigr)^{\einhalb}  }_{L^p(\calO)},  
  \]
  with a uniform constant.
\end{enumerate}

\section{Applications and consequences}

\subsection{Riesz sequences of sectorial exponentials}
\label{sec:riesz-exponentials}

A classical problem in harmonic analysis is to characterize those
sequences $(\lambda_n)_n$ in the right half-plane for which the
normalized exponentials
\[
  f_n(t)=\sqrt{2\Re\lambda_n}\,e^{-\lambda_n t}
\]
form a Riesz sequence in $L^2(\RR_+)$.  By means of the Laplace
transform, this problem is equivalent to the corresponding problem for
normalized reproducing kernels in the Hardy space of the right
half-plane.  The classical characterization
\cite{Carleson,ShapiroShields,NikolskiPavlov} is given by the uniform
separation condition
\[
  \inf_n\prod_{k\ne n} \left|\frac{\lambda_n-\lambda_k}    {\lambda_n+\overline{\lambda_k}}\right|>0.
\]
The connection with the Hautus condition has been exploited in the
work of Russell and Weiss \cite{RussellWeiss:Hautus} and Jacob and
Zwart
\cite{JacobZwart:diagonal-Hautus-Nikolskii-Pavlov,JacobZwart:SCL2001}.
As an application of our main result on Hilbert spaces, we inverse the
argument and give an alternative proof of the sectorial case of this
characterization which does not rely on interpolation theory in Hardy
spaces.  For positive real exponents, the same approach can be based
on the recent result of Daun and Jacob \cite[Theorem~1.5]{DaunJacob}
on self-adjoint generators.

\begin{proposition}\label{prop:riesz-exponentials}
  Let $\theta\in(0,\pihalbe)$ and let $(\lambda_n)_n$ be a sequence in
  $\sector{\theta}$.  Then the following assertions are equivalent:
  \begin{enumerate}
  \item The sequence $(\sqrt{2 \Re\lambda_n}\,e^{-\lambda_n\cdot})_n$ is
    a Riesz sequence in $L^2(\RR_+)$.
  \item The sequence $(\lambda_n)_n$ is separated in the
    pseudohyperbolic metric, that is,
    \begin{equation}\label{eq:exp-separation}
      \delta :=  \inf_{n\ne k}     \left|\frac{\lambda_n-\lambda_k}     {\lambda_n+\overline{\lambda_k}}\right|>0.
      \end{equation}
  \end{enumerate}
\end{proposition}

\begin{proof}
  The necessity of \eqref{eq:exp-separation} is classical and provided
  for convenience. Consider the normalised functions
  $f_n(t) = \sqrt{2 \Re(\lambda_n)}\, e^{-\lambda_n t}$ and observe
  that
\[
  d_{nk}^2 := 1-|[f_n,f_k]_{L^2}|^2
  =\frac{|\lambda_n+\overline{\lambda_k}|^2-4\Re\lambda_n\Re\lambda_k}{|\lambda_n+\overline{\lambda_k}|^2}
  =\left|\frac{\lambda_n-\lambda_k}{\lambda_n+\overline{\lambda_k}}\right|^2.
\]
The lower Riesz estimate now implies
$\norm{ \alpha f_n+\beta f_k }_{L^2}^2 \ge c(|\alpha|^2+|\beta|^2)$
which we can phrase in matrix language as
\[
G_{nk}=
\begin{pmatrix}
1 &[f_n,f_k]\\
\overline{[f_n,f_k]}&1
\end{pmatrix}
\ge c\, \mathrm{I}.
\]
The two eigenvalues of $G_{nk}$ are $1\pm |[f_n,f_k]|$, so that  $1-|[f_n,f_k]|\ge c$.
Inserting this into the above identity gives
\[
  d_{nk}^2 =1-|[f_n,f_k]|^2 \ge 1-(1-c)^2 =c(2-c) =: \delta^2 >0.
\]
  Conversely, assume \eqref{eq:exp-separation} and let
  $a_n=\Re\lambda_n$.  We consider the diagonal operator $A$ on
  $X=\ell^2$ given by $Ae_n=\lambda_n e_n$ and the observation
  operator
  \[
    Cx=\sum_n\sqrt{2 a_n}\,x_n.
  \]
  Clearly, $A$ is sectorial of angle at most $\theta$.  Moreover, as a
  diagonal operator, $A$ has a bounded $H^\infty$ functional calculus,
  and satisfies in particular lower square function estimates.  We
  first observe that
  \begin{equation}\label{eq:dyadic-count}
    \sup_{r>0}\quad \#\{n:r\le a_n<2r\}<\infty.
  \end{equation}
  Indeed, by sectoriality, the points $\lambda_n$
  with $r\le a_n<2r$ are contained in a rectangle
  of diameter comparable to $r$, whereas
  \eqref{eq:exp-separation} implies
  \[
    |\lambda_n-\lambda_k| \ge\delta(a_n+a_k)\ge 2\delta r
  \]
  for any two distinct such points.  A simple packing argument yields
  \eqref{eq:dyadic-count}. It follows that
  \[
    \sum_n \frac{a_n}{1+|\lambda_n|^2} = \sum_{k\in \ZZ} \sum_{2^k \le a_n < 2^{k+1} }  \frac{a_n}{1+|\lambda_n|^2} \lesssim  \sum_k \frac{2^k}{1+4^{k}}  <\infty,
  \]
  so that $C$ is bounded with respect to the graph norm of   $\DOMAIN(A)$.

  \medskip \noindent\emph{Admissibility.}  For fixed $n$ and
  $j\in\ZZ$, we let $ I_j(n) = \{ k:2^ja_n\le a_k<2^{j+1}a_n \}$.
  By \eqref{eq:dyadic-count}, the cardinalities of
  $I_j(n)$ are uniformly bounded by some constant $M$.
  Hence
  \[
    \sum_k \left|\frac{\sqrt{a_na_k}} {\lambda_n+\overline{\lambda_k}}\right|
    \le \sum_{j\in\ZZ}\sum_{k\in I_j(n)} \min\left\{\sqrt{\frac{a_n}{a_k}}, \sqrt{\frac{a_k}{a_n}}\right\}
   \le M\sum_{j\ge0}2^{-j/2}  +M\sum_{j<0}2^{(j+1)/2}  \lesssim M.
  \]
  By symmetry, Schur's test yields the admissibility of $C$.

  \medskip \noindent\emph{The Hautus condition.}
  Fix  $\lambda\in\CC_+$ and write $a=\Re\lambda$.  By
  \eqref{eq:exp-separation}, we may choose $\varepsilon\in(0,\einhalb)$,
  independently of $\lambda$, such that the disc $B(\lambda,\varepsilon a)$
  contains at most one point $\lambda_k$.  If the disc contains no
  such point, we immediately obtain
  \[
    \norm{ x }^2\le \frac{1}{\varepsilon^2a^2}\, \bignorm{ (\lambda-A)x    }^2.
  \]
  Otherwise, let $\lambda_k$ be the unique point in this disc. Then
  \begin{equation}\label{eq:exp-hautus-tail}
    \sum_{n\ne k}|x_n|^2    \le\frac{1}{\varepsilon^2a^2}    \norm{ (\lambda-A)x }^2.
  \end{equation}
  We claim that
  \begin{equation}\label{eq:exp-resolvent-sum}
    \sum_{n\ne k}\frac{a_n}{|\lambda-\lambda_n|^2}    \lesssim\frac1a.
  \end{equation}
  Indeed, splitting the sum according to whether $a_n<\frac{a}2$,  $\frac{a}2\le a_n<2a$, or $a_n\ge 2a$, we obtain
  \[
      \sum_{n\ne k}\frac{a_n}{|\lambda-\lambda_n|^2}
      \; \le  \;
      \Bigl(\frac4{a^2}\sum_{a_n<a/2}a_n\Bigr)    + \Bigl(\frac1{\varepsilon^2a^2}   \sum_{a/2\le a_n<2a,  n\ne k} a_n \Bigr)  +   \Bigl(4 \sum_{a_n\ge 2a}\frac1{a_n}\Bigr)
       \; \le   \;  \frac{4M}{a}(2+\varepsilon^{-2}).
  \]
  where we use \eqref{eq:dyadic-count} as well as dyadic counting for
  the first and last term, while the middle term contains at most $2M$
  terms of size $\sim a$.
  Writing
  \[
    \sqrt{a_k}\,x_k    =\frac1{\sqrt2}Cx    -\sum_{n\ne k}      \frac{\sqrt{a_n}}{|\lambda-\lambda_n|}      \cdot|\lambda-\lambda_n|x_n
  \]
  and using the Cauchy-Schwarz inequality together with
  \eqref{eq:exp-resolvent-sum}, we infer that
  \[
    \sqrt{a_k}|x_k|    \lesssim |Cx|+\frac1{\sqrt a}\norm{ (\lambda-A)x }.
  \]
  Since $a_k\simeq a$, this gives
  \[
    |x_k|\lesssim \frac1{\sqrt a}|Cx| +\frac1a\norm{ (\lambda-A)x },
  \]
  and combining  with \eqref{eq:exp-hautus-tail}, we arrive at the
  Hautus estimate
  \[
    \norm{ x }\lesssim \frac1{\Re\lambda}\norm{ (\lambda-A)x }
    +\frac1{\sqrt{\Re\lambda}}|Cx|.
  \]
  By Theorem~\ref{thm:Hautus-Hilbert}, $C$ is exactly observable,
  which is precisely the desired lower Riesz estimate.
\end{proof}

%
%
%
%
%
%
%
%
%
%

\appendix

\section{$\gamma$ functions. }

\begin{lemma}[Weis]\label{lem:weis}
  Let $X$ be a Banach space and $f: [a,b] \to X$ be differentiable on
  $(a,b)$ and have a Bochner integrable derivative.  Then
  $f \in \gamma(a,b; X)$ and
  \begin{equation}    \label{eq:weis-lemma}
    \norm{ f }_{\gamma(a, b; X)} \le  \sqrt{b-a} \Bigl( \norm{f(b) }_X + \int_a^b \norm{f'(s) }\,\mathrm{d}s\Bigr)
  \end{equation}
\end{lemma}
\begin{proof}
  We write $f(t) = b(b) - \int_t^b f'(s)\,\mathrm{d}s$. Since constant
  functions are radonifying, we only look at
  \begin{align*}
    \Bigl(\EE \Bignorm{ \sum_{n=N}^M \gamma_n \int_a^b f'(s) \int_a^s      e_n(t)\,\mathrm{d}t \, \mathrm{d}s }^2\Bigr)^\einhalb
    \le & \; 
    \int_a^b \norm{f'(s)}    \Bigl( \EE \Bigl| \sum_{n=N}^M \gamma_n [e_n, \eins_{[a,s]} ]    \Bigr|^2 \Bigr)^\einhalb \,\mathrm{d}s  \\
    = & \; 
    \int_a^b \norm{f'(s)} \Bigl(\sum_{n=N}^M | [e_n, \eins_{[a,s]} ]|^2 \Bigr)^\einhalb \,\mathrm{d}s.
  \end{align*}
  It follows that $\sum_n \gamma_n \int_a^b f(t)\,e_n(t)\,\mathrm{d}t$
  converges in $L^2(\Omega; X)$ as well as  \eqref{eq:weis-lemma}, as claimed.
\end{proof}

\begin{coro}\label{coro:W11}
  Let $X$ be a Banach space and $f \in W^{1,1}(\RR; X)$. Then $f$ is radonifying and
  \[
    \norm{f}_{\gamma(\RR; X) } \le 2 \norm{ f}_{W^{1,1}(\RR; X)}
  \]
\end{coro}
\begin{proof}
  Write $f(n+1) = f(t) + \int_t^{n+1} f'(s)\,\mathrm{d}s$ take norms and integrate in $t$ over $[n, n+1]$: it follows that
  \[
    \norm{ f(n+1) }_X \le \int_n^{n+1} \norm{ f(t) }_X \,\mathrm{d}t + \int_n^{n+1} \int_t^{n+1} \norm{ f'(s) }_X \,\mathrm{d}s \,\mathrm{d}t \le \int_n^{n+1} \norm{ f(t) }_X \,\mathrm{d}t + \int_n^{n+1}  \norm{ f'(s) }_X \,\mathrm{d}s
  \]
  Now, summing in $n$ shows that $\sum_n \norm{ f(n) }_X \le \norm{ f }_{W^{1,1}(\RR; X)}$ and the result follows from the preceding lemma.
\end{proof}

\def\SUBMITTED{Submitted}
\def\TOAPPEAR{To appear in }
\def\PREPARATION{In preparation }

\def\cprime{$'$}
\providecommand{\bysame}{\leavevmode\hbox to3em{\hrulefill}\thinspace}

\end{document}